\documentclass{amsart}

\usepackage{color}
\usepackage{epsfig}
\usepackage{amsthm}
\usepackage{amssymb}
\usepackage{amsmath}
\usepackage{amscd}

\newtheorem {Theorem}  {Theorem}

\numberwithin{Theorem}{section}

\newtheorem {Lemma}[Theorem]  {Lemma}
\newtheorem {Proposition}[Theorem]{Proposition}
\theoremstyle{definition}
\newtheorem{Definition}[Theorem]{Definition}
\theoremstyle{remark}
\newtheorem{Remark}[Theorem]{Remark}
\newtheorem{Example}[Theorem]{Example}

\newtheorem {Corollary}[Theorem]{Corollary}

\expandafter\chardef\csname pre amssym.def
at\endcsname=\the\catcode`\@ \catcode`\@=11
\def\undefine#1{\let#1\undefined}
\def\newsymbol#1#2#3#4#5{\let\next@\relax
 \ifnum#2=\@ne\let\next@\msafam@\else
 \ifnum#2=\tw@\let\next@\msbfam@\fi\fi
 \mathchardef#1="#3\next@#4#5}
\def\mathhexbox@#1#2#3{\relax
 \ifmmode\mathpalette{}{\m@th\mathchar"#1#2#3}%
 \else\leavevmode\hbox{$\m@th\mathchar"#1#2#3$}\fi}
\def\hexnumber@#1{\ifcase#1 0\or 1\or 2\or 3\or 4\or 5\or 6\or 7\or 8\or
 9\or A\or B\or C\or D\or E\or F\fi}

\font\teneufm=eufm10 \font\seveneufm=eufm7 \font\fiveeufm=eufm5
\newfam\eufmfam
\textfont\eufmfam=\teneufm \scriptfont\eufmfam=\seveneufm
\scriptscriptfont\eufmfam=\fiveeufm

\catcode`\@=\csname pre amssym.def at\endcsname

\newcounter{remark}
\def  \eps  {\epsilon}

\def\3n{\negthinspace \negthinspace \negthinspace }
\def\2n{\negthinspace \negthinspace }
\def\1n{\negthinspace }

\newcommand{\bg}{\begin{equation}}
\newcommand{\ed}{\end{equation}}
\newcommand{\bga}{\begin{eqnarray}}
\newcommand{\eda}{\end{eqnarray}}

\def\cbdu{\hfill{$\Box$}}

\renewcommand{\a}{\alpha}

\newcommand{\R}{\mathbf{R}}

\newcommand{\Mm}{{\mathcal M}}

\newcommand{\ueps}{u ^{\eps}}

\def\bex{\begin{Example}}
\def\eex{\end{Example}}

\def  \R   {{\mathbb R}}

\def  \12  {{\frac{1}{2}}}

\def\bd{\begin{Definition}}
\def\ede{\end{Definition}}
\def\be{\begin{equation}}
\def\bel{\begin{equation}\label}
\def\ee{\end{equation}}
\def\bt{\begin{Theorem}}
\def\et{\end{Theorem}}
\def\bc{\begin{Corollary}}
\def\ec{\end{Corollary}}
\def\bl{\begin{Lemma}}
\def\el{\end{Lemma}}
\def\bp{\begin{Proposition}}
\def\ep{\end{Proposition}}

\def\br{\begin{Remark}}
\def\er{\end{Remark}}
\def\ba{\begin{array}}
\def\ea{\end{array}}

\def\bea{\begin{eqnarray}}
\def\eea{\end{eqnarray}}

\begin{document}

\title[Decay of solutions for  approximations to the NSE]{Comparison of decay of solutions to two compressible approximations to Navier-Stokes equations}

\author{C\'esar J. Niche}
\address[C.J. Niche]{Departamento de Matem\'atica Aplicada, Instituto de Matem\'atica, Universidade Federal do Rio de Janeiro, CEP 21941-909, Rio de Janeiro - RJ, Brasil}
\email{cniche@im.ufrj.br}

\author{Mar\'{\i}a E. Schonbek}
\address[M.E. Schonbek]{Department of Mathematics, UC Santa Cruz, Santa Cruz, CA 95064, USA}
\email{schonbek@ucsc.edu}

\thanks{C.J. Niche acknowledges financial support from PRONEX  E-26/110.560/2010-APQ1, FAPERJ-CNPq and Ci\^encia sem Fronteiras - PVE 011/12. M. E. Schonbek was partially supported by NSF Grant DMS-0900909.}
\begin{abstract}
In this article, we use the decay character of initial data to compare the  energy decay rates of solutions to different compressible approximations to the Navier-Stokes equations. We show that the system having a nonlinear damping term has slower decay than its counterpart with an advection-like term. Moreover, me characterize a set of initial data for which the decay of the first system is driven by the difference between the full solution and the solution to the linear part, while for the second system the linear part provides the decay rate.
\end{abstract}

\date{\today}

\maketitle

\section{Introduction}

The Navier-Stokes equations for an incompressible homogeneous fluid in $\R^3$ 

\begin{eqnarray}
\label{eqn:navier-stokes}
\partial _t u + (u \cdot \nabla) u & = & \Delta u - \nabla p, \nonumber \\ div \, u & = & 0, \nonumber \\ u_0 (x) & = & u(x,0)
\end{eqnarray}
have been extensively studied because of their importance in modelling a wide range of phenomena in Fluid Mechanics. Taking divergence in the first line of (\ref{eqn:navier-stokes}), using the divergence-free condition and then  inverting the Laplacian, we obtain the nonlocal relation $p = - \Delta ^{-1} div \, (u \cdot \nabla) u$, which poses very hard problems when trying to solve these equations numerically. 

In order to avoid these problems, Temam \cite{MR0237972} proposed a model approximating (\ref{eqn:navier-stokes})  in which the pressure and the velocity are related through $\eps \, p = - div \, u$, for $\eps > 0$, thus breaking the nonlocality.  To ``stabilize'' this system, i.e. to have an energy inequality, he added the nonlinear term $\frac{1}{2} (div \, \ueps) \ueps$, which then leads to the compressible system

\begin{eqnarray}  \label{eqn:rusin}
\partial _t \ueps  + (\ueps \cdot \nabla) \ueps  + \frac{1}{2} (div \, \ueps) \ueps & = & \Delta \ueps + \frac{1}{\eps} \, \nabla \cdot div \, \ueps  \nonumber \\ \ueps (x,0)  & = & \ueps _0 (x).
\end{eqnarray}
This system has been used in many numerical experiments and has also been the subject of some articles concerning its analytical properties (see  Fabrie and Galusinski \cite{MR1868354}, Plech{\'a}{\v{c}} and {\v{S}}ver{\'a}k \cite{MR2012858}). Recently Rusin \cite{rusin} proved existence of global weak solutions in $\R^3$ and their convergence in $L^3 _{loc} (\R^3 \times \R_+) $, when $\eps$ goes to zero, to a suitable (in the sense of Caffarelli-Kohn-Nirenberg \cite{MR673830}) solution to the Navier-Stokes equation.

Given the differences in the linear parts and nonlinear terms of (\ref{eqn:navier-stokes}) and (\ref{eqn:rusin}) it is natural to ask whether these affect the decay rates of the $L^2$ norm of solutions. A useful tool to try to answer this question is the {\em decay character} of the initial datum, introduced by Bjorland and M.E. Schonbek \cite{MR2493562} and refined by Niche and M.E. Schonbek \cite{2014arXiv1405.7565N}. Roughly speaking, the decay character is a number associated to every $u_0 \in L^2 (\R^3)$ that describes the behaviour of $\widehat{u_0}$ near $\xi = 0$, which characterizes the norm decay of solutions to linear systems $u_t = {\mathcal L} u$ for a wide class of linear operators that includes $\mathcal{L} = \Delta$ and $\mathcal{L} = \Delta + \frac{1}{\eps}  \, \nabla \cdot div$. Using the sharp decay rate estimates obtained by Niche and M.E. Schonbek and the Fourier Splitting method, the following Theorem can be proved.

\bt{(Niche and M.E. Schonbek \cite{2014arXiv1405.7565N})} \label{decay-known-result} Let $u$ be a solution to either (\ref{eqn:navier-stokes}) or (\ref{eqn:rusin}), with  $u_0 \in L^2 (\R^3)$ and decay character $r^{\ast} = r^{\ast} (u_0)$, with $- \frac{3}{2} < r^{\ast} < \infty$. Then

\begin{displaymath}
\Vert u(t) \Vert _{L^2 (\R^3)} ^2 \leq C (1 + t) ^{- \min \{ \frac{3}{2} + r^{\ast}, \frac{5}{2} \}}.
\end{displaymath}
\et

\br \, Bjorland and M.E. Schonbek \cite{MR2493562} proved this  result for the Navier-Stokes equations (\ref{eqn:navier-stokes}).
\er

\smallskip 

The decay rate obtained in Theorem \ref{decay-known-result} provides plenty of information about the similarities between (\ref{eqn:navier-stokes}) and (\ref{eqn:rusin}). From Theorem \ref{characterization-decay-l2} we have that  

\begin{displaymath} 
\Vert e^{t {\mathcal L}} u_0 \Vert ^2 _{L^2 (\R^3)} \leq C (1 + t) ^{- \left( \frac{3}{2} + r^{\ast} \right)}, \qquad r^{\ast} = r^{\ast} (u_0),
\end{displaymath}
so we see that for $r^{\ast} \leq 1$ the linear parts are the ones that have slower decay, while for $r^{\ast} > 1$ the nonlinear terms are the ones driving  the decay due to the fast dissipation provided by $e^{t {\mathcal L}}$. Also, the stabilizing term $\frac{1}{2} (div \, \ueps) \ueps$, needed for having an energy inequality in (\ref{eqn:rusin}), does not change the relative strength of the linear parts and the nonlinear terms with regards to their influence on the decay of energy. We can then conclude that, regarding energy dissipation, the Navier-Stokes equations (\ref{eqn:navier-stokes}) and Temam's approximation (\ref{eqn:rusin}) have the same behaviour.

The main goal of this article is to study the decay of a different compressible approximation to (\ref{eqn:navier-stokes}) and compare its decay rates to those from Theorem \ref{decay-known-result} using the decay character. More precisely, consider the system

\begin{eqnarray} \label{eqn:lelievre-lemarieu} 
\partial _t \ueps  + (\ueps \cdot \nabla) \ueps  + \a \left| \ueps \right|^2 \ueps & = & \Delta \ueps + \frac{1} \eps \, \nabla \cdot div \, \ueps  \nonumber, \quad \a > 0 \\ \ueps (x,0)  & = & \ueps _0 (x),
\end{eqnarray}
introduced by Leli\`evre \cite{MR2833362} and Lemari\'e-Rieusset and Leli\`evre \cite{MR2835875} as a modification of one used by Vishik and Fursikov \cite{MR0609735} to construct statistical solutions to Navier-Stokes equations. System (\ref{eqn:lelievre-lemarieu}) has the same scaling as Navier-Stokes equations, thus allowing for its analysis in spaces which contain homogeneous initial data leading to selfsimilar solutions. Note that this system differs from Temam's approximation (\ref{eqn:rusin}) in that, instead of having an advection-like term $(div \, \ueps) \ueps$, it has a strongly nonlinear damping term $\left| \ueps \right|^2 \ueps$. Leli\`evre \cite{MR2833362} proved that for $\ueps _0 \in L^2 (\R^3)$, there exists a global in time weak solution to (\ref{eqn:lelievre-lemarieu})  and that when $\a$ and $\eps$ go to zero, solutions to (\ref{eqn:lelievre-lemarieu}) converge, as distributions, to a suitable solution to the Navier-Stokes equations.

We now state our main result which will allow us to compare the behaviour of (\ref{eqn:navier-stokes}) and (\ref{eqn:rusin}) to that of (\ref{eqn:lelievre-lemarieu}).

\bt
\label{main-theorem}
Let $\ueps _0 \in L^2 (\R^3)$, with $r^{\ast} = r^{\ast} (\ueps _0)$, with $- \frac{3}{2} < r^{\ast} < \infty$. Then for any weak solution to (\ref{eqn:lelievre-lemarieu}) we have that for $\a > \frac{\eps}{4}$, $C = C(\a, \eps, \Vert u_0 ^{\eps} \Vert _{L^2})$ and any $\delta > 0$

\begin{displaymath}
\Vert \ueps (t) \Vert ^2 _{L^2} \leq C (1 + t) ^{- \min \{ \frac{3}{2} + r^{\ast}, \frac{3}{2} - \delta \} }.
\end{displaymath}
\et

Thus, the strong nonlinear damping $\left| \ueps \right|^2 \ueps$ also leads to decay (see the energy inequality (\ref{eqn:energy-inequality-llr})). However, {\em there are significant quantitative and qualitative differences regarding the decay of solutions}. First, note that for initial data with $r^{\ast} \geq 0$, solutions to Navier-Stokes equations (\ref{eqn:navier-stokes}) and Temam's approximation (\ref{eqn:rusin})  have faster decay rates than solutions to Lemari\'e-Rieusset and Leli\`evre's system. This is so precisely due to the presence of this nonlinear damping, which slows down decay, as can be explicitly seen in the proof by comparing (\ref{eqn:estimate-advection-term}) to (\ref{eqn:estimate-damping}).

Moreover, for  Navier-Stokes equations and Temam's approximation, the linear part determines the decay rates for $r^{\ast} < 1$, while for Lemari\'e-Rieusset and Leli\`evre's system the linear part is the leading one only for $r^{\ast} < 0$. Hence,  for  $\ueps _0$ with $0 < r^{\ast} < 1$ not only the decay rates are different, but so is  the dissipation mechanism, given by the linear part in the first case, and by the nonlinear terms in the second case. In Remark \ref{special-data} we show that $v_0 = \Lambda^s \ueps _0$, where $\ueps _0$ is in $L^1 (\R^3)  \cap H^s (\R^3)$ with $0 < s < 1$, has decay character $r^{\ast} (v_0) = s$, thus providing explicit examples of initial data that lead to such behaviour.

\br For any fixed $\a > 0$, when $\eps$ goes to zero we obtain the Navier-Stokes equations with an extra damping term. The decay of solutions to this system has been recently addressed, see Cai and Lei \cite{Xiaojing20101235}, Jia, Zhang and Dong \cite{Jia20111736}, Jiang and Zhu \cite{MR2871810}, Jiang \cite{MR2927562}. Theorem \ref{main-theorem} improves and generalizes the results obtained in the case of $\beta = 3$ in \cite{MR2927562}.
\er

This article is organized as follows. In Section \ref{settings} we recall some existence results and properties of solutions to (\ref{eqn:lelievre-lemarieu}) and, following Niche and M.E. Schonbek  \cite{2014arXiv1405.7565N}, we provide the definitions and results we need concerning the decay character and the characterization of decay of linear systems. In Section \ref{proof}, we prove our main Theorem \ref{main-theorem}.

\section{Settings}

\label{settings}

\subsection{Solutions to (\ref{eqn:lelievre-lemarieu})} For the sake of completeness, we first recall results concerning existence of solutions to (\ref{eqn:lelievre-lemarieu}) and their convergence to solutions to (\ref{eqn:navier-stokes}).

\bt{(Theorem 3.3, Leli\`evre \cite{MR2833362})} \label{thm-existence-lelievre} Let $u_0 \in L^2 (\R^3)$. Then there exists a distributional solution $\ueps$ to (\ref{eqn:lelievre-lemarieu}) such that $\ueps \in L^{\infty} (\R_+, L^2) \cap L^2 (\R_+, \dot{H}^1) \cap L^4 (\R_+, L^4)$.
\et

\bt{(Theorem 4.3, Leli\`evre \cite{MR2833362})} Let $u_0 \in L^2 (\R^3)$, with $\nabla \cdot u_0 = 0$. Then solutions from Theorem  \ref{thm-existence-lelievre} converge when $\a$ and $\eps$ go to zero, as distributions, to solutions to (\ref{eqn:navier-stokes}).
\et 

\subsection{Decay character} As the long time behaviour of solutions to many dissipative systems is determined by the low frequencies of the solution, Bjorland and M.E. Schonbek \cite{MR2493562} introduced the idea of {\em decay character} of a function $u_0$ in $L^2 (\R^n)$ in order to characterize the decay of solutions to Navier-Stokes equations with that initial datum. Recently, Niche and M.E. Schonbek \cite{2014arXiv1405.7565N} generalized this notion in order to use data in $H^s (\R^n), s  > 0$ and to obtain results for other equations. We recall now these definitions and results.

\bd Let  $u_0 \in L^2(\R^n)$ and $\Lambda = (- \Delta) ^{\frac{1}{2}}$. For $s \geq 0$ the {\em s-decay indicator}  $P^s _r (u_0)$ corresponding to $\Lambda ^s u_0$ is 

\begin{displaymath}
P^s _r (u_0) = \lim _{\rho \to 0} \rho ^{-2r-n} \int _{B(\rho)} |\xi|^{2s}| \widehat{u}_0 (\xi)| ^2 \, d \xi
\end{displaymath}
for $r \in \left(- \frac{n}{2} + s, \infty \right)$, where $B(\rho)$ is the ball at the origin with radius $\rho$.
\ede

\br \label{heuristics} \, Setting $r = q + s$, we see that the $s$-decay indicator compares $|\widehat{\Lambda ^s u_0} (\xi)|^2$ to $f(\xi) = |\xi|^{2(q + s)}$ near $\xi = 0$. When $s = 0$ we recover the definition of decay indicator given by Bjorland and M.E. Schonbek \cite{MR2493562}.
\er

\bd  \label{df-decay-character} The {\em decay character of $\Lambda^s u_0$}, denoted by $r_s ^{\ast} = r_s ^{\ast}( u_0)$ is the unique  $r \in \left( -\frac{n}{2} + s, \infty \right)$ such that $0 < P^s _r (u_0) < \infty$, provided that this number exists. If such  $P^s _r ( u_0)$ does not exist, we set $r_s ^{\ast} = - \frac{n}{2} + s$, when $P^s _r (u_0)  = \infty$ for all $r \in \left( - \frac{n}{2} + s, \infty \right)$  or $r_s ^{\ast} = \infty$, if $P^s _r (u_0)  = 0$ for all $r \in \left( -\frac{n}{2} + s, \infty \right)$. 
\ede

\br \label{remark-lp-l2} \, (Examples 2.5 and 2.6, Niche and M.E. Schonbek \cite{2014arXiv1405.7565N}).  Let $u_0 \in L^2 (\R^n)$ such that $\widehat{u}_0 (\xi) = 0$, for $|\xi| < \delta$, for some $\delta > 0$. Then, $P_r ^s (u_0) = 0$, for any $r \in \left( - \frac{n}{2} + s, \infty\right)$ and  $r^{\ast} _s (u_0) = \infty$.  If $u_0 \in L^p (\R^n) \cap L^2 (\R^n)$, with $1  \leq p \leq 2$, then $r^{\ast} (u_0) = - n \left( 1 - \frac{1}{p} \right)$, so if $u_0 \in L^1 (\R^n) \cap L^2 (\R^n)$ we have that $r^*(u_0) =0$ and if $u_0 \in L^2 (\R^n)$ but $u_0 \notin L^p (\R^n)$, for any $1 \leq p < 2$, we have that $r^*(u_0) = -  \frac{n}{2}$.
\er

\smallskip

Let $u_0 \in H^s (\R^n)$, with $s > 0$. As $\widehat{\Lambda^s (u_0)} (\xi) = |\xi|^s \widehat{u_0} (\xi)$, the heuristics for the decay character given in Remark \ref{heuristics} lead us to expect that $r^{\ast} _s (u_0) = r^{\ast} (\Lambda ^s u_0) = s + r^{\ast} (u_0)$. This is the content of the following Theorem.

\bt \label{decay-character-hs} (Theorem 2.11, Niche and M.E. Schonbek \cite{2014arXiv1405.7565N}) Let $u_0 \in H^s (\R^n), s > 0$.  
\begin{enumerate}
 \item If $-\frac{n}{2} < r^{\ast} (u_0) < \infty$ then  $- \frac{n}{2} +s< r_s^{\ast}(u_0) < \infty$ and  $r_s^{\ast}(u_0) = s + r^{\ast} (u_0)$. 
\item $r_s^{\ast}  (u_0) = \infty$ if and only if $r^{\ast} (u_0) = \infty$.
\item  $r^{\ast} (u_0) =- \frac{n}{2}$ if and only if $r_s^{\ast}(u_0) = r^{\ast} (u_0) + s = - \frac{n}{2} + s$. \end{enumerate}
\et

\br \label{special-data} We can now justify the assertion made after the statement of Theorem \ref{main-theorem} about initial data for which the decay of solutions to (\ref{eqn:rusin}) and  (\ref{eqn:lelievre-lemarieu})  are quantitatively and qualitatively different. Let $u_0 \in L^1 (\R^n) \cap L^2 (\R^n)$ such that $\Lambda ^s u_0 \in L^2 (\R^n)$, with $0 < s < 1$. From Remark \ref{remark-lp-l2} we have that $r^{\ast} (u_0) = 0$, while from Theorem  \ref{decay-character-hs} we have that $r^{\ast} (\Lambda ^s u_0) = r^{\ast} _s (u_0) =  s + r^{\ast} (u_0) = s$. This proves that $\Lambda^s u_0$ has decay character $0 < r^{\ast} < 1$.
\er

\subsection{Linear operators and characterization of decay} We describe now the linear operators for which we characterize the decay of the $L^2$ norm of solutions in terms of the decay character. For a Hilbert space $X$ on $\R^n$, we consider a pseudodifferential operator $\mathcal{L}: X^n \to \left( L^2 (\R^n) \right) ^n$, with symbol $ \Mm(\xi)$ such that 

\be
\label{eqn:symbol}
\Mm(\xi) = P^{-1} (\xi) D(\xi) P(\xi), \qquad \xi-a.e.
\ee
where $P(\xi) \in O(n)$ and $D(\xi) = - c_i |\xi|^{2\a} \delta _{ij}$, for $c_i > c>0$ and $0 < \a \leq 1$.  Taking the Fourier Transform of the linear equation

\be
\label{eqn:linear-part}
v_t = \mathcal{L} v,
\ee
multiplying by $\widehat{v}$, integrating in space and then using (\ref{eqn:symbol}) we obtain

\begin{displaymath}
\frac{1}{2} \frac{d}{dt} \Vert v(t) \Vert _{L^2} ^2 \leq  - C  \int _{\R^n} |\xi|^{2 \a} |\widehat{v}|^2 \, d \xi
\end{displaymath}
which is the key inequality for using the Fourier Splitting method in the proofs. 

\br \, The fractional Laplacian on vector fields in $\R^n$

\begin{displaymath}
\left( \mathcal{L} u \right) _i = (- \Delta) ^{\a} u_i, \qquad i = 1, \cdots n,
\end{displaymath}
provides an example of (\ref{eqn:symbol}), as its symbol is $\left( \Mm(\xi) \right) _{ij}= - C |\xi|^{2 \a} \delta _{ij}$. The operator 

\begin{displaymath}
\mathcal{L} u = \Delta u + \frac{1}{\eps} \, \nabla div \, u , \qquad \eps > 0,
\end{displaymath}
i.e. the linear part of (\ref{eqn:lelievre-lemarieu}), provides a second example, as (\ref{eqn:symbol}) holds with $ \left( \Mm(\xi) \right) _{ij} = - |\xi|^2 \delta_{ij}  - \frac{1}{\eps} \xi_i \xi _j $, $D(\xi) = diag (- |\xi|^2, - |\xi|^2, - \left( 1 + \frac{1}{\eps} \right) |\xi|^2 )$ and

\begin{displaymath}
P(\xi) =  \left( \begin{array} {ccc}
\frac{- \xi_2}{\sqrt{\xi_1 ^2 + \xi_2 ^2}} & \frac{- \xi_1 \xi_3}{\sqrt{1 - \xi_3 ^2}} & \xi_1 \\ \frac{\xi_1}{\sqrt{\xi_1 ^2 + \xi_2 ^2}} &  \frac{- \xi_2 \xi_3}{\sqrt{1 - \xi_3 ^2}} & \xi_2 \\ 0 &  \frac{1 - \xi_3 ^2}{\sqrt{1 - \xi_3 ^2}} & \xi_3
\end{array} \right),
\end{displaymath}
where $v = (\xi_1, \xi_2, \xi_3)$ has norm one. As a result of this, 

\be
\label{eqn:sol-fund-frec}
\left(  e^{t \Mm (\xi)} \right) _{ij} = e ^{-t |\xi|^2} \delta _{ij} - \frac{\xi _i \xi _j}{|\xi|^2} \left(e ^{-t |\xi|^2} - e ^{- \left( 1 + \frac{1}{\eps} \right)t |\xi|^2} \right),
\ee
see Rusin \cite{rusin}.
\er

\smallskip

We now state the Theorem that describes decay in terms of the decay character for linear operators as in (\ref{eqn:symbol}).

\bt{(Theorem 2.10, Niche and M.E. Schonbek \cite{2014arXiv1405.7565N})}
\label{characterization-decay-l2}
Let $v_0 \in L^2 (\R^n)$ have decay character $r^{\ast} (v_0) = r^{\ast}$. Let $v (t)$ be a solution to  (\ref{eqn:linear-part}) with data $v_0$. Then:
\begin{enumerate}
\item if $- \frac{n}{2 } < r^{\ast}< \infty$, there exist constants $C_1, C_2> 0$ such that
\begin{displaymath}
C_1 (1 + t)^{- \frac{1}{\a} \left( \frac{n}{2} + r^{\ast} \right)} \leq \Vert v(t) \Vert _{L^2} ^2 \leq C_2 (1 + t)^{- \frac{1}{\a} \left( \frac{n}{2} + r^{\ast} \right)};
\end{displaymath}
\item if $ r^{\ast}= - \frac{n}{2}$, there exists $C = C(\eps) > 0$ such that
\begin{displaymath}
\Vert v(t) \Vert _{L^2} ^2 \geq C (1 + t)^{-\eps}, \qquad \forall \eps > 0,
\end{displaymath}
i.e. the decay of $\Vert v(t) \Vert _{L^2} ^2$ is slower than any uniform  algebraic rate;
\item if $r^{\ast} = \infty$, there exists $C > 0$ such that
\begin{displaymath}
\Vert v(t) \Vert _{L^2} ^2 \leq C (1 + t) ^{- m}, \qquad \forall m > 0,
\end{displaymath}
i.e. the decay of $\Vert v(t) \Vert _{L^2}$ is faster than any algebraic rate.
\end{enumerate} 
\et

\section{Proof of Theorem \ref{main-theorem}}

\label{proof}

\noindent
{\bf Proof:} \,  The proof is based on the Fourier Splitting method, introduced by M.E. Schonbek to study decay of parabolic conservation laws \cite{MR571048} and of Navier-Stokes equations \cite{MR775190}, \cite{MR837929}. As is usual in this context, we prove the estimate assuming the solutions are regular enough (for the existence of regular approximations, see Leli\`evre \cite{MR2833362}). The limiting argument that  proves the estimate for weak solutions follows that for the Navier-Stokes equations in pages 267--269 in Lemari\'e-Rieusset \cite{MR1938147} and the Appendix in Wiegner \cite{MR881519}, which we refer to for full details.

We first show that solutions obey an energy inequality. As

\begin{displaymath}
\int _{\R^3} \ueps \left( \ueps \cdot \nabla \right) \ueps \, dx = - \frac{1}{2} \int _{\R^3} |\ueps|^2 div \, \ueps \, dx,
\end{displaymath}
it follows that

\begin{eqnarray}
\label{eqn:estimate-energy-inequality}
\left| \int _{\R^3} \ueps \left( \ueps \cdot \nabla \right) \ueps \, dx \right| & \leq & \frac{1}{2} \left| \int _{\R^3} |\ueps|^2 div \, \ueps \, dx \right| = \frac{1}{2} \Vert |\ueps|^2 div \, \ueps \Vert _{L^1} \nonumber \\ & \leq & \frac{1}{2} \Vert |\ueps|^2 \Vert _{L^2}  \Vert div \, \ueps \Vert _{L^2} \leq \frac{\eps}{4}  \Vert \ueps \Vert _{L^4} ^4 +  \frac{1}{4 \eps} \Vert div \, \ueps \Vert _{L^2} ^2 .
\end{eqnarray}
Multipyling (\ref{eqn:lelievre-lemarieu}) by $\ueps$, integrating in space and using (\ref{eqn:estimate-energy-inequality}) we obtain 

\be
\label{eqn:energy-inequality-llr}
\frac{1}{2} \frac{d}{dt} \Vert \ueps (t) \Vert _{L^2} ^2 \leq - \Vert \nabla \ueps (t) \Vert _{L^2} ^2 -  \frac{3}{4 \eps}  \Vert div \, \ueps \Vert _{L^2} ^2 - \left( \a - \frac{\eps}{4} \right) \Vert \ueps (t) \Vert _{L^4} ^4,
\ee
which, given the hypotheses on $\eps$ and $\a$, provides an energy inequality. Then 

\be
\label{eqn:ei}
\frac{1}{2} \frac{d}{dt} \Vert \ueps (t) \Vert _{L^2} ^2 \leq - \Vert \nabla \ueps (t) \Vert _{L^2} ^2 -  C \Vert div \, \ueps \Vert _{L^2} ^2 .
\ee
Now let

\begin{displaymath}
B(t) = \{\xi \in \R^3: |\xi|^2 \leq \frac{r'(t)}{2 C r(t)} \}
\end{displaymath}
where $r$ is a positive increasing function with $r(0) = 1$ and $C$ is an appropiate constant. From (\ref{eqn:ei}) we have that

\begin{displaymath}
\frac{1}{2} \frac{d}{dt} \Vert \ueps (t) \Vert _{L^2} ^2 \leq - C \int _{\R^3} |\xi|^2 |\widehat{\ueps} (\xi, t)|^2 \, d \xi \leq - \frac{r'(t)}{2 r(t)} \int _{B(t) ^c} |\widehat{\ueps} (\xi, t)|^2 \, d \xi . 
\end{displaymath}
Adding and substracting a term similar to the one on the right side of this inequality, only that with $B(t)$ as the domain of integration, and then mutiplying by $r(t)$ we obtain

\begin{equation}
\label{eqn:fourier-splitting-main}
\frac{d}{dt} \left( r(t) \Vert \ueps (t) \Vert _{L^2} \right) \leq r'(t) \int _{B(t)}   |\widehat{\ueps} (\xi, t)|^2 \, d \xi.
\end{equation}
We now prove a pointwise estimate for 

\begin{displaymath}
\widehat{\ueps} (\xi, t) = e^{t \mathcal{M} (\xi)} \widehat{\ueps _0} (\xi)- \int _0 ^t e^{(t - s) \mathcal{M} (\xi)} G(\xi, s) \, ds
\end{displaymath}
where $e^{t \mathcal{M} (\xi)}$ is as in (\ref{eqn:sol-fund-frec}) and 
\begin{displaymath}
G(\xi, s) = {\mathcal F} \left( (\ueps \cdot \nabla) \ueps + \left| \ueps \right|^2 \ueps  \right),
\end{displaymath}
where $\mathcal{F}$ is also the Fourier transform. As 

\begin{displaymath}
(\ueps \cdot \nabla) \ueps = \nabla \cdot (\ueps \otimes \ueps) - (div \, \ueps)  \ueps
\end{displaymath}
and 

\begin{displaymath}
\left| \widehat{(div \, \ueps)  \ueps} (\xi, t)\right| \leq \Vert div \, \ueps (t) \Vert _{L^2}  \Vert \ueps(t) \Vert _{L^2},
\end{displaymath}
we obtain

\begin{displaymath}
|{\mathcal F} \left( (\ueps \cdot \nabla) \ueps (t) \right) | \leq |\xi| \Vert \ueps(t) \Vert _{L^2} ^2 + \frac{1}{2} \Vert div \, \ueps (t) \Vert _{L^2}  \Vert \ueps (t) \Vert _{L^2} \leq C |\xi| \Vert \ueps (t) \Vert _{L^2} ^2.
\end{displaymath}
Also

\begin{displaymath}
\left| {\mathcal F} \left(  \left| \ueps \right|^2 \ueps   \right) (t) \right| \leq \Vert {\mathcal F} \left(  \left| \ueps \right|^2 \ueps   \right) (t) \Vert _{L^{\infty}} \leq \Vert \left| \ueps (t)  \right|^2 \ueps (t)  \Vert _{L^1} \leq \Vert \ueps (t) \Vert _{L^3} ^3.
\end{displaymath}
Now, we estimate the nonlinear transport term. We have
\begin{eqnarray}
\label{eqn:inequality-integral-term}
\left| \int _0 ^t e^{(t -s) \mathcal{M} (\xi)} {\mathcal F} \left( (\ueps \cdot \nabla) \ueps \right) (\xi, s) \, ds \right| & \leq & C \int _0 ^t e^{- C (t - s) |\xi|^2} \, |\xi| \Vert \ueps(s) \Vert _{L^2} ^2 \, ds \nonumber \\ & \leq & C |\xi| \left( \int _0 ^t \Vert \ueps(s) \Vert _{L^2} ^2 \, ds \right).
\end{eqnarray} 
Suppose now that 

\be
\label{eqn:estimate-beta}
\Vert \ueps(t) \Vert _{L^2} ^2 \leq C (1 + t) ^{- \beta},
\ee
for  some $\beta \geq 0$ and $\beta \neq 1$. We then have, after choosing $r(t) = (t + 1)^3$ 

\begin{eqnarray}
\label{eqn:estimate-advection-term}
\int _{B(t)} \left( \int _0 ^t e^{(t -s) \mathcal{M} (\xi)} {\mathcal F} \left( (\ueps \cdot \nabla) \ueps \right) (\xi, s) \, ds \right) ^2 d \xi & \leq & C \left( \frac{r'(t)}{r(t)} \right) ^{\frac{5}{2}} (1 + t) ^{2 (1 - \beta)} \nonumber \\ & \leq & C (1 + t) ^{- \left( \frac{1}{2} +2 \beta \right)},
\end{eqnarray}
where we used (\ref{eqn:inequality-integral-term}) and (\ref{eqn:estimate-beta}). We now estimate the damping term. We have 

\begin{align}
\label{eqn:estimate-damping}
\int _{B(t)} \left( \int _0 ^t e^{(t -s) \mathcal{M} (\xi)} {\mathcal F} \left(  \left| \ueps \right|^2 \ueps  \right) (\xi, s) \, ds \right) ^2 d \xi  \leq C \int _{B(t)} \left( \int _0 ^t \Vert \ueps (s) \Vert _{L^3} ^3 \, ds \right) ^2 \, d \xi \nonumber \\  \leq  C \int _{B(t)} \left( \int _0 ^t \Vert \ueps (s) \Vert_{L^2}  \Vert \ueps (s) \Vert _{L^4} ^2 \, ds \right) ^2 \, d \xi \nonumber \\  \leq  C \int _{B(t)} \left( \int _0 ^t \Vert \ueps (s) \Vert_{L^2} ^2 \, ds \right) \left( \int _0 ^t  \Vert \ueps (s) \Vert _{L^4} ^4 \, ds \right)  d \xi \nonumber \\  \leq  C (t + 1) ^{- \left( \frac{1}{2} + \beta \right)} 
\end{align}
where we used the interpolation
\begin{displaymath}
\Vert \ueps (t) \Vert _{L^3}  \leq \Vert \ueps (t) \Vert_{L^2} ^{\frac{1}{3}} \Vert \ueps (t)  \Vert _{L^4} ^{\frac{2}{3}},
\end{displaymath}
H\"older's inequality, (\ref{eqn:estimate-beta}) and the fact that $\ueps \in L^{\infty} _t L^2 _x \cap L^4 _t L^4 _x$. Of these two terms,  (\ref{eqn:estimate-damping}) has the slower decay. 

The only apriori estimate we have is $\Vert \ueps (t) \Vert _{L^2} \leq C$, i.e. $\beta = 0$. So, in the ball $B(t)$ we have 

\begin{eqnarray}
\int _{B(t)} |\widehat{\ueps} (\xi, t)| ^2 \, d \xi & \leq & C \int _{B(t)} |e^{t \mathcal{M} (\xi)} \widehat{u_0} |^2 \, d \xi \nonumber \\ & + & C \int _{B(t)} \left( \int _0 ^t e^{(t - s) \mathcal{M} (\xi)} G(\xi, s) \, ds \right)^2 \, d \xi \nonumber \\
& \leq & C (t + 1) ^{- \left( \frac{3}{2} + r^{\ast} \right)} + C (t + 1) ^{- \frac{1}{2}}, 
\end{eqnarray}
which leads to 

\be
\Vert \ueps (t) \Vert _{L^2} ^2 \leq C (1 + t) ^{- \min \{ \frac{3}{2} + r^{\ast}, \frac{1}{2} \}}.
\ee
We now use this to boostrap our decay estimates by improving the value of $\beta$. When $r^{\ast} < -1$, we have that $\beta = \frac{3}{2} + r^{\ast}$ and (\ref{eqn:estimate-damping}) provides no improvement. When $r^{\ast} \geq -1$, then $\beta = \frac{1}{2}$ which leads to 

\begin{eqnarray}
\int _{B(t)} |\widehat{\ueps} (\xi, t)| ^2 \, d \xi & \leq & C \int _{B(t)} |e^{t \mathcal{M} (\xi)} \widehat{u_0} |^2 \, d \xi \nonumber \\ & + & C \int _{B(t)} \left( \int _0 ^t e^{(t - s) \mathcal{M} (\xi)} G(\xi, s) \, ds \right)^2 \, d \xi \nonumber \\
& \leq & C (t + 1) ^{- \left( \frac{3}{2} + r^{\ast} \right)} + C (t + 1) ^{- 1}.
\end{eqnarray}
Hence

\be
\Vert \ueps (t) \Vert _{L^2} ^2 \leq C (1 + t) ^{- \min \{ \frac{3}{2} + r^{\ast}, 1 \}}.
\ee
Again, $r^{\ast} < - \frac{1}{2}$ does not produce an improvement on the decay while $r^{\ast} \geq - \frac{1}{2}$ leads to $\beta = 1$, for which (\ref{eqn:estimate-advection-term}) and (\ref{eqn:estimate-damping}) do not hold, as

\begin{displaymath}
\int _0 ^t \Vert \ueps (s) \Vert _{L^2} ^2 \, ds \leq C \ln (t + 1).
\end{displaymath}
However, working along the lines of these estimates we obtain

\begin{equation}
\int _{B(t)} \left( \int _0 ^t e^{(t -s) \mathcal{M} (\xi)} {\mathcal F} \left( (\ueps \cdot \nabla) \ueps \right) (\xi, s) \, ds \right) ^2 d \xi  \leq  C \left(  \frac{r'(t)}{r(t)} \right) ^{\frac{5}{2}} \ln ^2 (t + 1)
\end{equation}
and

\begin{equation}
\int _{B(t)} \left( \int _0 ^t e^{(t -s) \mathcal{M} (\xi)} {\mathcal F} \left(  \left| \ueps \right|^2 \ueps  \right) (\xi, s) \, ds \right) ^2 d \xi  \leq C \left(  \frac{r'(t)}{r(t)} \right) ^\frac{3}{2} \ln (t + 1).
\end{equation}
Again, the second term has slower decay and choosing $r(t) = (t + 1) ^{\gamma}$ for large enough $\gamma > 0$ we have that, as $\ln (t + 1) \leq (t + 1) ^{\delta}$ for any $\delta > 0$  

\begin{eqnarray}
\int _{B(t)} |\widehat{\ueps} (\xi, t)| ^2 \, d \xi & \leq & C \int _{B(t)} |e^{t \mathcal{M} (\xi)} \widehat{u_0} |^2 \, d \xi \nonumber \\ & + & C \int _{B(t)} \left( \int _0 ^t e^{(t - s) \mathcal{M} (\xi)} G(\xi, s) \, ds \right)^2 \, d \xi \nonumber \\
& \leq & C (t + 1) ^{- \left( \frac{3}{2} + r^{\ast} \right)} + C (t + 1) ^{- \left( \frac{3}{2} - \delta \right)}.
\end{eqnarray}
Then

\begin{displaymath}
\Vert \ueps (t) \Vert _{L^2} ^2 \leq C (1 + t) ^{- \min \{ \frac{3}{2} + r^{\ast}, \frac{3}{2} - \delta \}}, \qquad \forall  \delta > 0.
\end{displaymath}
This boostrapping process does not improve the decay further, as now for $r^{\ast} \geq 0$, the integral of $\Vert \ueps (t) \Vert _{L^2}$ is bounded by a constant. $\Box$ 

\bibliographystyle{unsrt}
\bibliography{StronglyDamped.bib}

\end{document}